\documentclass[12pt]{amsart}
\usepackage{amssymb,amsmath,color}
\usepackage{graphicx,amsmath}
\usepackage{listings}
\usepackage{pst-node}
\usepackage{tikz-cd} 
\usepackage{hyperref}
\usepackage[nameinlink]{cleveref}

\usepackage{mathrsfs}

\newcommand{\tn}[1]{\textnormal{#1}}
\newcommand{\R}[0]{\mathbb{R}}
\newcommand{\C}[0]{\mathbb{C}}

\newcommand{\Z}[0]{\mathbb{Z}}
\newcommand{\<}[0]{\langle}
\renewcommand{\>}[0]{\rangle}
\newcommand{\func}[3]{#1:#2\rightarrow #3}
\newcommand{\colvec}[2][1.0]{%
  \scalebox{#1}{%
    \renewcommand{\arraystretch}{1.0}%
    $\begin{bmatrix}#2\end{bmatrix}$%
  }
}
\makeatletter
\renewcommand*\env@matrix[1][*\c@MaxMatrixCols c]{%
  \hskip -\arraycolsep
  \let\@ifnextchar\new@ifnextchar
  \array{#1}}
\makeatother

\usepackage{soul}

\theoremstyle{plain}
\newtheorem{thm}{Theorem}[section]
\newtheorem{lemma}[thm]{Lemma}
\newtheorem*{theorem*}{Theorem}
\newtheorem*{corollary*}{Corollary}

\newtheorem{prop}[thm]{Proposition}

\theoremstyle{definition}

\newtheorem{remark}[thm]{Remark}

\numberwithin{equation}{section}

\usepackage{vmargin}
\setpapersize{USletter}

\setmargrb{2.5cm}{1.5cm}{2.5cm}{2cm}

 \usepackage{scalerel,stackengine}
\stackMath
\newcommand\reallywidehat[1]{%
\savestack{\tmpbox}{\stretchto{%
  \scaleto{%
    \scalerel*[\widthof{\ensuremath{#1}}]{\kern.1pt\mathchar"0362\kern.1pt}%
    {\rule{0ex}{\textheight}}
  }{\textheight}%
}{2.4ex}}%
\stackon[-6.9pt]{#1}{\tmpbox}%
}
\title{Distribution of eigenvalues of the Kohn Laplacian on Sphere Quotients}

\author{Adam Cohen}
\address[Adam Cohen]{School of Mathematics, University of Minnesota, Minneapolis, MN 55455, USA}
\email{cohen795@umn.edu}

\author{Yash Rastogi}
\address[Yash Rastogi]{School of Mathematics, Georgia Institute of Technology, Atlanta, GA 30332, USA}
\email{yrastogi7@gatech.edu}

\author{Samuel Sottile}
\address[Samuel Sottile]{Department of Mathematics, Stanford University, Stanford, CA 94305, USA}
\email{sottile@stanford.edu}

\author{Yunus E. Zeytuncu}
\address[Yunus E. Zeytuncu]{Department of Mathematics and Statistics, University of Michigan--Dearborn, Dearborn, MI 48128, USA}
\email{zeytuncu@umich.edu}

\thanks{This work is supported by NSF (DMS-1950102 and DMS-2243808) and NSA (H98230-23).}
\begin{document}

\nocite{*}
\begin{abstract}
We prove a Weyl-type theorem for the Kohn Laplacian on sphere quotients as CR manifolds. We show that we can determine the fundamental group from the spectrum of the Kohn Laplacian in dimension three. Furthermore, we prove Sobolev estimates for the complex Green's operator on these quotient manifolds.
\end{abstract}

\maketitle

\section{Introduction}
\subsection{Motivation}
A well-known question in Riemannian geometry, popularized in Mark Kac's article ``Can One Hear the Shape of a Drum?'' \cite{kacdrum}, is whether two Riemannian manifolds with the same Dirichlet eigenvalues are necessarily isometric. Milnor found counterexamples in dimension $16$ \cite{Milnor1964} and later Gordon and Webb found counterexamples for planar domains \cite{GordonWebb}. Ikeda produced counterexamples that fail to be homotopy equivalent: five dimensional lens spaces \cite{IYCounterexample}. While extensively studied for the Laplace-Beltrami operator, the isospectral drum question is also of interest beyond the setting of Riemannian manifolds. Recent research has addressed this question in the context of CR geometry. Folland provides crucial tools for understanding the relationship between isospectrality and isometry on sphere quotients \cite{FollandSpheres}, leveraging the representation theory of the unitary group to compute the eigenvalues and eigenfunctions of the Kohn Laplacian $\Box_b$ on $(0,q)$-forms. This approach facilitates the proof of asymptotic estimates for the distribution of eigenvalues of the Kohn Laplacian on spheres \cite{REU2020WeylSpheres} and enables the construction of counterexamples akin to those in \cite{IYCounterexample}, as demonstrated in \cite{GLR}.

Weyl's law is a celebrated result in Riemannian geometry that relates the distribution of the Dirichlet eigenvalues of the Laplace-Beltrami operator on a bounded domain to that domain's volume. Obtaining a comparable result for CR manifolds is nontrivial since the Kohn Laplacian is not elliptic. Stanton and Tartakoff established a version for $(p,q)$-forms with $0 < q < n$ for $(2n+1)$-dimensional strongly pseudoconvex CR manifolds \cite{ST84} and recent work has established a version for functions on lens spaces \cite{REU2021Lens}. We extend these results by producing a Weyl-type theorem for the Kohn Laplacian on functions on all sphere quotients. This demonstrates that we can discern the order of the fundamental groups of sphere quotients, and thereby ``hear'' their volume. We further investigate the spectrum in dimension $3$ as there are infinitely many pairs of $(2n-1)$-dimensional CR lens spaces which are CR isospectral but not CR isometric for $n \geq 3$ \cite{GLR}. We show that one can hear the fundamental group for three-dimensional sphere quotients.

Given $f \in (\ker \Box_b)^{\perp}$, a weak solution to the Poisson equation for the Kohn Laplacian $\Box_b u = f$ is given by $u = \mathcal{G} f + g$, where $\mathcal{G}$ is the complex Green's operator and $g \in \ker ( \Box_b )$. The equation $\Box_b u = f$ governs the behavior of boundary values of holomorphic functions and is discussed further in \cite{Chen-Shaw}. Folland and Kohn show $\mathcal{G}$ gives one weak derivative for certain $q$ values based on the CR structure \cite{FollandKohn}. Through spectral analysis, we offer an elementary proof of this fact for functions on sphere quotients. The reader should note that this is not covered by the estimate in \cite{FollandKohn}. 

In \Cref{sec:weyl}, we prove a Weyl-type theorem for the Kohn Laplacian on sphere quotients, thus demonstrating that we can hear the order of the fundamental group. In \Cref{sec:dim3}, we demonstrate that we can hear the fundamental group of spherical $3$-manifolds. We extend the work of \cite{REU2019Sobolev} to provide Sobolev estimates for the complex Green's operator on sphere quotients in \Cref{sec:Sobolev}. 

\subsection{Background}
The unitary group acts faithfully on $S^{2n-1}$ by left multiplication. With the induced metric and CR structure from $S^{2n-1} \subset \C^n$, the action is by CR isometries. This descends to a left action on $L^2(S^{2n-1})$ via precomposition with the inverse map. We have the decomposition 
\[L^2(S^{2n-1}) = \bigoplus_{p,q}\mathcal{H}_{p,q}(S^{2n-1}),\]
where $\mathcal{H}_{p,q}(S^{2n-1})$ denotes the space of harmonic homogeneous polynomials of bidegree $(p,q)$ \cite{FollandSpheres}. For simplicity, we denote $\mathcal{H}_{p,q}(S^{2n-1})$ as $\mathcal{H}_{p,q}$ when the dimension $n$ is clear from the context. Each $\mathcal{H}_{p,q}$ is an irreducible representation of $U(n)$ and is an eigenspace for $\Box_b$ \cite{FollandSpheres}.

Let $g \in U(n)$. Since $g$ is a normal matrix, by the spectral theorem, $g$ is conjugate by unitary matrices to the diagonal matrix
\[\mu = \colvec{\mu_1 & \cdots & 0\\ \vdots & \ddots & \vdots \\ 0 & \cdots & \mu_n}\]
with eigenvalues $\mu_1,\ldots, \mu_n$. Consequently, for the character $\chi_{\mathcal{H}_{p,q}}$, we have
\[\chi_{\mathcal{H}_{p,q}}(\mu) = \chi_{\mathcal{H}_{p,q}}(g).\]
With an appropriate choice of basis for $\mathcal{H}_{p,q}$, computing this character is straightforward. One such basis is given by
\[\{f_{\alpha,\beta} = \overline{D}^\alpha D^\beta \|z\|^{2-2n}, \alpha,\beta \tn{ multiindices, } |\alpha|=p,|\beta|=q, \alpha_1=0\tn{ or }\beta_1=0\}\]
where
\[\overline{D}^\alpha = \frac{\partial^{\alpha_1}}{\partial \overline{z_1}^{\alpha_1}}\frac{\partial^{\alpha_2}}{\partial \overline{z_2}^{\alpha_2}}\cdots \frac{\partial^{\alpha_n}}{\partial \overline{z_n}^{\alpha_n}} = \frac{\partial^{|\alpha|}}{\partial \overline{z}_1^{\alpha_1}\cdots\partial\overline{z}_n^{\alpha_n}},\qquad D^\beta = \frac{\partial^{|\beta|}}{\partial z_1^{\beta_1}\cdots\partial z_n^{\beta_n}}.\]
For more details, we refer the reader to \cite{AxlerHFT}, \cite{REU2017Rossi}, and \cite{REU2021Lens}. We call the $(\alpha,\beta)$ that meet the requirements above admissible for $(p,q)$. By the chain rule, we compute
\begin{align*}
    \frac{\partial}{\partial z_i}\|z\|^{2-2n}|_{z=w} = \frac{\partial}{\partial z_i}\|\mu^{-1} z\|^{2-2n}|_{z=w} =\mu_i^{-1}\frac{\partial}{\partial z_i}\|z\|^{2-2n}|_{z=\mu^{-1}w}.
\end{align*}
Writing $\mu^\alpha = \mu_1^{\alpha_1} \mu_2^{\alpha_2}\cdots \mu_n^{\alpha_n}$ and similarly for $\beta$, we have
\[\overline{\mu}^\alpha\mu^\beta f_{\alpha,\beta} = f_{\alpha,\beta}\circ \mu^{-1}.\] 
Therefore, 
\[\chi_{\mathcal{H}_{p,q}}(g) = \chi_{\mathcal{H}_{p,q}}(\mu) = \sum_{\alpha,\beta}\overline{\mu}^\alpha\mu^\beta,\]
where the sum is taken over all admissible $(\alpha,\beta)$ pairs. Let $G \subset U(n)$ be a discrete group acting freely on $S^{2n-1}$. We define $\chi_{\mathcal{H}_{p,q}}^G$ to be the character restricted to $G$. Let $\func{\pi}{S^{2n-1}}{G \backslash S^{2n-1}}$ be the quotient map. One can compute
\[\dim \pi^* L^2(G \backslash S^{2n-1}) \cap \mathcal{H}_{p,q} = \< \chi_{\mathcal{H}_{p,q}}^G,\chi^G_{\tn{Triv}}\>_G = \frac{1}{|G|}\sum_{g \in G}\chi_{\mathcal{H}_{p,q}}(g).\]
 We will additionally denote this subspace of $\mathcal{H}_{p,q}$ fixed by the action of $G$ as $\mathcal{H}_{p,q}^G$. Note that $\pi^*$ is injective and commutes with $\Box_b$, so the spectrum is given by the $\mathcal{H}_{p,q}^G$. We note that if $G,G' \subset U(n)$ are two finite subgroups that are conjugate, then their respective quotients are CR isometric.

We compute the characters of the $\mathcal{H}_{p,q}$ spaces for $S^3$ both to illustrate the general case and for use in \Cref{sec:dim3}. Let $g \in U(2)$ and let $\mu_1,\mu_2$ be the eigenvalues of $g$. Then 
\[\chi_{\mathcal{H}_{p,q}}(g) = \sum_{\alpha,\beta} \overline{\mu_1}^{\alpha_1}\overline{\mu_2}^{\alpha_2}\mu_1^{\beta_1}\mu_2^{\beta_2} = \sum_{\alpha,\beta}\mu_1^{\beta_1-\alpha_1}\mu_2^{\beta_2-\alpha_2}\]
The admissible $(\alpha, \beta)$ pairs in the sum can be enumerated as follows: 
\[\left(\colvec{p\\0},\colvec{0\\q} \right),\left(\colvec{p-1\\1},\colvec{q\\0}\right),\ldots , \left(\colvec{0\\p},\colvec{0\\q}\right), \left(\colvec{0\\p},\colvec{1\\q-1}\right),\ldots , \left(\colvec{0\\p},\colvec{q\\0}\right).\]
So, 
\[\chi_{\mathcal{H}_{p,q}}(g) = \frac{\mu_2^q}{\mu_1^p}\sum_{j=0}^{p+q}\left( \frac{\mu_1}{\mu_2} \right)^j.\]
If $g \in SU(2)$, then $\mu_2 = \mu_1^{-1}$. So we have
\[\chi_{\mathcal{H}_{p,q}}(g) = \mu_1^{-p-q}\sum_{j=0}^{p+q}\mu_1^{2j}.\]

The finite subgroups of $SU(2)$ can be classified up to conjugation into two infinite families and three exceptional groups. Since $SU(2) \cong S^3$, these groups act freely on $S^3$. Considering the double cover $S^3\rightarrow SO(3)$, these arise from the finite subgroups of $SO(3)$. The finite subgroups of $SO(3)$ are the cyclic groups $\Z/m\Z$, the dihedral groups $D_{2m}$, and three exceptional groups, the tetrahedral group $T \cong A_4$, the octahedral group $O \cong S_4$, and the icosahedral group $I \cong A_5$. Under the double cover $SU(2) \rightarrow SO(3)$, these groups lift to finite groups in $SU(2)$ with twice the order. These groups together with the odd order cyclic groups are all of the finite subgroups of $SU(2)$. The lifts of the dihedral groups are called the binary dihedral groups $2D_{2m}$. One should note $2D_{2} \cong \Z/4\Z$. We will assume $m \geq 2$ when talking about the groups $2D_{2m}$. The lifts of the exceptional groups are $2T$, $2O$, and $2I$, the binary tetrahedral, octahedral, and icosahedral groups respectively. The quotients $(\Z/m\Z)\backslash S^3$ are the lens spaces $L(m;1,-1)$. The case of quotients of the sphere by a cyclic group was considered in \cite{REU2021Lens}.

The eigenvalue of $\Box_b$ on $\mathcal{H}_{p,q}$ is $2q(p+n-1)$ \cite{FollandSpheres}. So the multiplicity of an eigenvalue $\lambda$ for $\Box_b$ on $L^2(G \backslash S^{2n-1})$ is
\[\sum\dim \mathcal{H}_{p,q}^G,\]
where the sum is over those $p$ and $q$ satisfying $\lambda = 2q(p+n-1)$. Sobolev and Schatten estimates for the complex Green's operator for $\Box_b$ on $L^2(S^{2n-1})$ are in \cite{REU2019Sobolev}.

\subsection{Main Results} In \Cref{sec:weyl}, we prove a Weyl-type theorem, which relates the asymptotics of the eigenvalue counting function for the Kohn Laplacian on a sphere quotient to the volume of the sphere quotient.

Let $G \subset U(n)$ be a discrete group acting freely on $S^{2n-1}$. Let $N_{G\backslash S^{2n-1}}(\lambda)$ be the number of positive eigenvalues (counting multiplicity) of $\Box_b$ on $L^2(G \backslash S^{2n-1})$ that are less than or equal to $\lambda$.
\begin{thm}\label{thm:weyl's law} There exists a constant $C$ such that for every finite subgroup of $G$ of $U(n)$,
    \[\lim_{\lambda \rightarrow \infty} \frac{N_{G\backslash S^{2n-1}}(\lambda)}{\lambda^n} = C \frac{\tn{Vol}(S^{2n-1})}{|G|} = C \, \tn{Vol} (G \backslash S^{2n-1}).\]
\end{thm}
\begin{remark}
    Note the following behavior of the limit in \Cref{thm:weyl's law}. When the metric on the sphere $S^{2n-1}$ is scaled by a factor of $\beta$, the limit is scaled by a factor of $\beta^{2n}$. The Kohn Laplacian only depends on the sub-Riemannian structure of $S^{2n-1}$, i.e. \[\left(S^{2n-1}, \tn{Re}(\mathbb{L} \oplus \overline{\mathbb{L}}), g|_{\tn{Re}(\mathbb{L} \oplus \overline{\mathbb{L}})}\right).\] This has Hausdorff dimension $2n$ (more generally, see \cite{MitchellSubriemannian}), which is consistent with the behavior of the limit under scaling of the metric.
\end{remark}

In \Cref{sec:dim3}, we show the following result.
\begin{thm}\label{thm:hear pi1}
    We can determine the isomorphism class of $G \subset U(2)$ from the spectrum of $\Box_b$ on $L^2(G\backslash S^3)$.
\end{thm}
This is attained through a series of calculations involving computing multiplicities of specific eigenvalues of the Kohn Laplacian on sphere quotients.
\begin{remark}
    If $G$ and $G'$ are non-abelian, isomorphic subgroups of $U(2)$ that act freely on $S^3$, then they are conjugate in $U(2)$. This is not true for abelian subgroups of $U(2)$. If one could use the spectrum of the Kohn Laplacian to distinguish three-dimensional lens spaces with the same order fundamental group, then CR isometry would imply CR isospectrality for three-dimensional sphere quotients.
\end{remark}
In \Cref{sec:Sobolev}, we show that the the complex Green's operator gives one weak derivative.

\section{Estimates for the Dirichlet Eigenvalues of the Kohn Laplacian}\label{sec:weyl}

Let $N_{S^{2n-1}}(\lambda)$ be the number of positive eigenvalues for $\Box_b$ on $S^{2n-1}$ that are less than or equal to $\lambda$.
The group $U(n)$ acts by CR isometries on the sphere $S^{2n-1}$. If a discrete group $G \subset U(n)$ acts freely on $S^{2n-1}$, then the quotient is a CR manifold. For $\lambda > 0$, define
\[\mathcal{H}(\lambda) = \bigoplus_{0 < q(p+n-1) \leq \lambda}\mathcal{H}_{p,q}.\]
For $G \subset U(n)$ a discrete subgroup, define
\[\mathcal{H}^G(\lambda) = \bigoplus_{0 < q(p+n-1) \leq \lambda} \mathcal{H}^G_{p,q}.\]
Then
\[N_{S^{2n-1}}(2\lambda) = \dim \mathcal{H}(\lambda).\]
\[N_{G \backslash S^{2n-1}}(2\lambda) = \dim \mathcal{H}^G(\lambda).\]
Note that $\mathcal{H}^G(\lambda)$ still makes sense when $G$ does not act freely. The goal of this section to prove the following lemma.
\begin{lemma}\label{lem:weyl 2023}
    Let $G \subset U(n)$ be a finite subgroup. Then
    \[\lim_{\lambda \rightarrow \infty} \frac{\dim\mathcal{H}(\lambda)}{\dim\mathcal{H}^G(\lambda)} = |G|.\]
\end{lemma}
\begin{thm}[\cite{REU2020WeylSpheres}]\label{thm:weyl 2020}
    \[\lim_{\lambda \rightarrow \infty} \frac{N_{S^{2n-1}}(\lambda)}{\lambda^n} = \tn{Vol}(S^{2n-1}) \frac{n-1}{n(2\pi)^n\Gamma (n+1)}\int_{-\infty}^\infty \left( \frac{\tau}{\sinh \tau} \right)^ne^{-(n-2)\tau}d\tau.\]
\end{thm}
\Cref{thm:weyl's law} follows immediately from \Cref{lem:weyl 2023} and \Cref{thm:weyl 2020}.
We define
\[\Xi_\lambda = \left[ \sum_{k=0}^{\lfloor \lambda \rfloor -n+1}\binom{k+n-2}{n-2}\left[ \binom{\left\lfloor \frac{\lambda}{k+n-1}\right\rfloor+n-2}{n-1}\right]+ \sum_{k=1}^{\left\lfloor \frac{\lambda}{n-1} \right\rfloor}\binom{k+n-2}{n-2}\binom{\left\lfloor \frac{\lambda}{k} \right\rfloor}{n-1}\right].\]
\begin{lemma}\label{lem:halfway weyl estimate}
    \[\left | N_{G \backslash S^{2n-1}}(2\lambda) - \frac 1 {|G|}N_{S^{2n-1}}(2\lambda) \right| \leq (|G| - 1)\Xi_\lambda.\]
\end{lemma}
This lemma will be proved at the end of this section. Once we have this lemma, to prove \Cref{lem:weyl 2023}, we need to show:
\[\lim_{\lambda \rightarrow \infty} \frac{\Xi_\lambda}{\dim \mathcal{H}(\lambda)} = 0.\]
\begin{proof}(Of \Cref{lem:weyl 2023})
First we observe from \Cref{thm:weyl 2020} that $\mathcal{H}(\lambda) \approx C\lambda^n$ as $\lambda \rightarrow \infty$. To show the above limit, we will show $\Xi_\lambda$ is $o(\lambda^n)$. Observe
\[\binom{a+b}{b} \leq (b+1)a^b\]
for $a \geq 1$. Also
\[\binom{a+b}{b} \leq (a+1)^b,\] \[\binom{a}{b} \leq a^b.\]
So
\[\Xi_\lambda \leq \left[ \sum_{k=0}^{\lfloor \lambda \rfloor - n+1}\binom{k+n-2}{n-2}\left(\frac{\lambda}{k+n-1} \right)^{n-1} + \sum_{k=1}^{\left\lfloor \frac{\lambda}{n-1} \right\rfloor}\binom{k+n-2}{n-2}\left(\frac{\lambda}{k} \right)^{n-1} \right]\]
\[\leq \lambda^{n-1} + \sum_{k=1}^{\lfloor \lambda \rfloor - n + 1}(n-1)k^{n-2}\left(\frac{\lambda}{k} \right)^{n-1} + \sum_{k=1}^{\left\lfloor \frac{\lambda}{n-1} \right\rfloor}(n-1)k^{n-2}\left( \frac{\lambda}{k} \right)^{n-1}.\]
We compute the inequality
\[\sum_{k=1}^x\frac{1}{k} \leq 1 + \int_{1}^{x-1}\frac 1x dx = 1 + \log(x-1).\]
So,
\[\Xi_\lambda \leq \lambda^{n-1} + (n-1)\lambda^{n-1}(1+ \log(\lambda - n)) + (n-1)\lambda^{n-1}\left(1+\log\left(\frac{\lambda}{n-1} -1 \right)\right)\]
Hence, $\Xi_\lambda$ is $O(\lambda^{n-1}\log \lambda)$.
\end{proof}

\begin{proof} (Of \Cref{lem:halfway weyl estimate})\\
For a multiindex $\alpha$, let $\tilde{\alpha} = (\alpha_1,\ldots, \alpha_{n-1},0)$.
Let $g \in U(n)$. Then $g$ has eigenvalues $\mu_1,\ldots, \mu_n$. We compute:
\begin{align*}
    \chi_{\mathcal{H}(\lambda)}(g) &= \sum_{\alpha_1 = 0,\ |\alpha| \leq \lfloor \lambda \rfloor - n+1}\mu^{-\alpha} \left[ \sum_{ |\beta| \leq \left\lfloor \frac{\lambda}{|\alpha| + n-1}  \right\rfloor,\ \beta_1 \geq 1}\mu^\beta  \right]\\
    &+\sum_{\beta_1=0,\ 1 \leq |\beta| \leq \left\lfloor \frac{\lambda}{n-1} \right\rfloor} \mu^\beta\left[ \sum_{|\alpha|\leq \left\lfloor \frac{\lambda}{|\beta|}\right\rfloor - n + 1} \mu^{-\alpha}\right].
\end{align*}
The first term counts the $f_{\alpha,\beta} = \overline{D}^\alpha D^\beta |z|^{2-2n}$ with $\alpha_1 = 0$, $\beta_1 \neq 0$, the second term counts those with $\beta_1=0$.

To compute $N_{G \backslash S^{2n-1}}$, we have:
\[N_{G \backslash S^{2n-1}}(2\lambda) = \dim \mathcal{H}^G(\lambda) = \frac{1}{|G|}\sum_{g \in G}\chi_{\mathcal{H}(\lambda)}(g) = \frac 1{|G|}\dim \mathcal{H}(\lambda) + \frac{1}{|G|}\sum_{g \in G, g \neq 1}\chi_{\mathcal{H}(\lambda)}(g).\]
Note that if $g$ is not the identity, then it must have an eigenvalue not equal to $1$. We can assume this is $\mu_n$.
We rewrite
\begin{align*}
    \chi_{\mathcal{H}(\lambda)}(g) &= \sum_{\alpha_1=0,\ |\alpha| \leq \lfloor \lambda \rfloor - n+1}\mu^{-\alpha} \left[\sum_{|\tilde{\beta}| \leq \left\lfloor \frac{\lambda}{|\alpha| + n-1}  \right\rfloor,\ \beta_1 \geq 1}\mu^{\tilde{\beta}}\sum_{\beta_n=0}^{\left\lfloor \frac{\lambda}{|\alpha| + n-1} \right\rfloor - \left|\tilde{\beta}\right|}\mu_n^{\beta_n}  \right]\\
    &+\sum_{\beta_1=0,\ 1 \leq |\beta| \leq \left\lfloor \frac{\lambda}{n-1} \right\rfloor} \mu^{\beta}\left[ \sum_{|\tilde{\alpha}| \leq \left\lfloor \frac{\lambda}{|\beta|}\right\rfloor - n + 1} \mu^{-\tilde{\alpha}}\sum_{\alpha_n=0}^{\left\lfloor \frac{\lambda}{|\beta|}\right\rfloor - n + 1 - |\tilde{\alpha}|}\mu_n^{-\alpha_n}\right].
\end{align*}
If $\zeta \neq 1$ is a $k$-th root of unity, then we have the trivial bound
\[\left | \sum_{j=a}^b \zeta^j\right| < k.\]
We know that $\mu_n \neq 1, \mu_n^{|G|} = 1$.
Using the above estimate with a liberal application of the triangle inequality and the `sticks and stones' combinatorial formula, we get
\begin{align*}
    |\chi_{\mathcal{H}(\lambda)}(g)|&\leq \sum_{k=0}^{\lfloor \lambda \rfloor -n+1}\binom{k+n-2}{n-2}\sum_{j=1}^{\left\lfloor \frac{\lambda}{k+n-1}\right\rfloor} \binom{j-1+n-2}{n-2}|G|\\
    &+\sum_{k=1}^{\left\lfloor \frac{\lambda}{n-1} \right\rfloor}\binom{k+n-2}{n-2}\sum_{j=0}^{\left\lfloor \frac{\lambda}{k} \right\rfloor-n+1}\binom{j + n-2}{n-2}|G|\\
    &= \sum_{k=0}^{\lfloor \lambda \rfloor -n+1}\binom{k+n-2}{n-2} \binom{\left\lfloor \frac{\lambda}{k+n-1}\right\rfloor+n-2}{n-1}|G|+ \sum_{k=1}^{\left\lfloor \frac{\lambda}{n-1} \right\rfloor}\binom{k+n-2}{n-2}\binom{\left\lfloor \frac{\lambda}{k} \right\rfloor}{n-1}|G|\\
    &= \Xi_\lambda |G|.
\end{align*}
We have
\[N_{G \backslash S^{2n-1}}(2\lambda) = \frac{1}{|G|}N_{G\backslash S^{2n-1}}(2\lambda) + \frac 1{|G|} \sum_{g \in G, g\neq 1}\chi_{\mathcal{H}(\lambda)}(g).\]
So we see \[\left | N_{G \backslash S^{2n-1}}(2\lambda) - \frac 1 {|G|}N_{S^{2n-1}}(2\lambda) \right| \leq \frac 1 {|G|}\sum_{g \in G, g \neq 1}|\chi_{\mathcal{H}(\lambda)}(g)| \leq (|G| - 1)\Xi_\lambda.\]
\end{proof}

\section{Classification of Spherical \texorpdfstring{$3$}{3}-Manifolds through Spectral Analysis}\label{sec:dim3}
In this section, we provide a proof by casework of \Cref{thm:hear pi1}. First, we give convenient presentations of the finite subgroups of $SU(2)$. Then we compute representatives for each of the conjugacy classes of the finite subgroups of $U(2)$ that act freely on the sphere, in a form conducive to computing their characters.

\subsection{CR Isometry for \texorpdfstring{$SU(2)$}{SU(2)} Quotients}\label{subsec:su2}
We show that in dimension 3, CR isospectral sphere quotients by finite subgroups of the special unitary group are CR isometric.
We use the identification of $SU(2)$ with $S^3$ considered as the unit quaternions. A unit quaternion $g = a+bi+cj+dk$ corresponds to
\[\colvec{a+bi & -c + di\\ c+di & a-bi}.\]
Then the trace of $g$ as a matrix is twice the real part of $g$ as a quaternion. This allows us to determine the eigenvalues of $g$ from the real part of $g$ since the determinant of $g$ is $1$. Recall that the finite subgroups of $SU(2)$ can be classified as the cyclic groups $\Z/m\Z$, the binary dihedral groups $2D_{2m}$, and the three exceptional groups, $2T,2O,2I$, the binary tetrahedral, octahedral and icosahedral groups respectively. For each of these groups, we pick a representative of each conjugacy class and systematically denote the representatives with the symbols $\Z/m\Z$, $2D_{2m}$, $2T,\ 2O$ and $2I$. Using convenient presentations of these groups, we compute their characters and subsequently compute the dimensions of the invariant $\mathcal{H}_{p,q}$ spaces.

The cyclic group $\Z/m\Z$ is given by the group generated by the quaternion $e^{2\pi i/m}$. The eigenvalues of this generator are $e^{\pm 2 \pi i/m}$. So we compute
\[\dim \mathcal{H}_{p,q}^{\Z/m\Z} = \frac 1m \sum_{n=0}^{m-1}e^{-(p+q)2\pi i n/m}\sum_{l=0}^{p+q}e^{4\pi i n l/m}.\]
For $m$ even, we compute:
\[\dim \mathcal{H}_{p,q}^{\Z/m\Z} = \begin{cases}
2 \left\lfloor \frac{p+q}{m} \right\rfloor + 1 & p+q \tn{ even}\\
0 & p+q \tn{ odd}
\end{cases}.\]
For $m$ odd, we compute:
\[\dim \mathcal{H}_{p,q}^{\Z/m\Z} = \begin{cases}
    2 \left\lfloor \frac{p+q}{2m} \right\rfloor + 1 & p+q \tn{ even}\\
 2 \left\lfloor \frac{p+q + m}{2m} \right\rfloor & p+q \tn{ odd}
\end{cases}.\]

The group $2D_{2m}$ is given by the set of quaternions
\[2D_{2m} = \{e^{\pi i l/m}, je^{\pi i l/m} \mid l = 0,\cdots 2m-1\}.\]
Since $2m$ of these elements have real part zero, their eigenvalues are $\pm i$. The other elements form the cyclic subgroup $\Z/2m\Z$. Then we compute
\[\tn{dim} \mathcal{H}_{p,q}^{2D_{2m}} = \frac 12 \tn{dim} \mathcal{H}_{p,q}^{\Z/2m\Z} + \frac 12 i^{-p-q}\sum_{l=0}^{p+q}(-1)^l = \frac 12 \tn{dim} \mathcal{H}_{p,q}^{\Z/2m\Z} + \frac 12\begin{cases}
    1 & p+q \equiv 0 \tn{ mod } 4\\
    -1 & p+q \equiv 2 \tn{ mod } 4\\
    0 & \tn{else}
\end{cases}. \]
We observe that when $p+q$ is odd, $\tn{dim} \mathcal{H}_{p,q}^{2D_{2m}} = 0$.

The binary tetrahedral group $2T$ has $24$ elements. It contains the quaternion group $Q = 2D_4$ and the remaining $16$ elements come from the vectors in $\R^4 \cong \mathbb{H}$ given by $\frac 12 (\epsilon_0, \epsilon_1,\epsilon_2,\epsilon_3)$, where $\epsilon_l \in \{\pm 1\}$. The eigenvalues of these $16$ elements not in $Q$ are $\frac 12\epsilon_0 \pm \frac{\sqrt{3}}{2}i$. Then we compute
\begin{align*}
    \dim \mathcal{H}_{p,q}^{2T} &= \frac 13 \dim \mathcal{H}_{p,q}^Q + \frac 13\left[e^{-\pi i (p+q)/3}\sum_{l=0}^{p+q}e^{2\pi i l/3} + e^{-2\pi i (p+q)/3}\sum_{l=0}^{p+q}e^{-2\pi i l/3}  \right]\\
&=\frac 13 \dim \mathcal{H}_{p,q}^Q + \frac 13\left[e^{\pi i (p+q)/3}\sum_{l=0}^{p+q}e^{2\pi i (l-(p+q))/3} + e^{-2\pi i (p+q)/3}\sum_{l=0}^{p+q}e^{-2\pi i l/3}  \right]\\
&=\frac 13 \dim \mathcal{H}_{p,q}^Q + \frac 13\left[(e^{\pi i (p+q)/3}+e^{-2\pi i (p+q)/3})\sum_{l=0}^{p+q}e^{-2\pi i l/3}  \right]\\
&= \frac 13 \dim \mathcal{H}_{p,q}^Q + \frac 13\begin{cases}
    2 & p+q \equiv 0 \tn{ mod } 6\\
    -2 & p+q \equiv 4 \tn{ mod } 6\\
    0 & \tn{else}
\end{cases}.
\end{align*}
The binary octahedral group $2O$ has $48$ elements. It contains the binary tetrahedral group $2T$ and the remaining $24$ elements are given by $\frac 1 {\sqrt{2}}(\epsilon_0, \epsilon_1, 0,0)$, and all of their permutations, where $\epsilon_l \in \{\pm 1\}$. Twelve of these elements have eigenvalues $\pm i$, six of these elements have eigenvalues $\frac 1{\sqrt{2}} \pm \frac 1{\sqrt{2}}i = e^{\pm \pi i/8}$, and the remaining six elements have eigenvalues $- \frac 1 {\sqrt{2}} \pm \frac 1{\sqrt{2}} i = e^{\pm 3\pi i/8}$. From this we can compute
\begin{align*}
    \dim \mathcal{H}^{2O}_{p,q} &= \frac 12 \dim\mathcal{H}^{2T}_{p,q}+\frac 14 i^{-p-q}\sum_{l=0}^{p+q}(-1)^l + \frac 18 \left[ e^{-\pi i(p+q)/8}\sum_{l=0}^{p+q}i^l + e^{-3\pi i(p+q)/8}\sum_{l=0}^{p+q}(-i)^l \right]\\
    &= \frac 12 \dim\mathcal{H}^{2T}_{p,q} + \frac 14\begin{cases}
    1 & p+q \equiv 0 \tn{ mod } 4\\
    -1 & p+q \equiv 2 \tn{ mod } 4\\
    0 & \tn{else}
\end{cases}
+\frac 18 \begin{cases}
    2 & p+q \equiv 0,2 \tn{ mod } 8\\
    -2 & p+q \equiv 4,6 \tn{ mod } 8\\
    0 & \tn{else}
\end{cases}\\
&= \frac 12 \dim\mathcal{H}^{2T}_{p,q}+\frac 12 \begin{cases}
    1 & p+q \equiv 0 \tn{ mod } 8\\
    -1 & p+q \equiv 6 \tn{ mod } 8\\
    0 & \tn{else}
\end{cases}.
\end{align*}
The binary icosahedral group $2I$ has $120$ elements. It contains the binary tetrahedral group $2T$ and the remaining $96$ elements are the even permutations of $\frac 12(0,\epsilon_1,\epsilon_2\phi^{-1},\epsilon_3\phi)$, where $\epsilon_l \in \{\pm 1\}$, and $\phi = \frac{1+\sqrt{5}}{2}$ is the golden ratio. The eigenvalues in $SU(2)$ are:
\begin{center}
\begin{tabular}{ |c|c| } 
 \hline
Real part of $g$ & Eigenvalues in $SU(2)$\\
\hline
$0$ & $\pm i$ \\ 
\hline
$\frac{1}{2}$ & $\frac{1}{2}\pm \frac{i\sqrt{3}}{2}$ \\ 
\hline
$-\frac{1}{2}$ &  $-\frac{1}{2}\pm \frac{i\sqrt{3}}{2}$ \\ 
\hline
$\frac{1}{2\phi}$ & $e^{\pm 2\pi i/5}$  \\ 
\hline
$-\frac{1}{2\phi}$ & $e^{\pm 3\pi i/5}$   \\ 
\hline
$\frac{\phi}{2}$ & $e^{\pm\pi i/5}$   \\ 
\hline
$-\frac{\phi}{2}$ & $e^{\pm 4\pi i/5}$   \\  
 \hline
\end{tabular}
\end{center}
From this, we compute:
\begin{align*}
    \dim \mathcal{H}_{p,q}^{2I} &= \frac 15 \dim \mathcal{H}_{p,q}^{2T} + \frac 15 i^{-p-q}\sum_{l=0}^{p+q}(-1)^k\\
    &+ \frac 1{10}(e^{\pi i (p+q)/3}+e^{-2\pi i (p+q)/3})\sum_{l=0}^{p+q}e^{-2\pi i l/3}+\frac 1{10}\sum_{m=1}^4e^{-\pi i l(p+q)/5}\sum_{l=0}^{p+q}e^{2\pi i m l/5}\\
  &= \frac 15 \dim \mathcal{H}_{p,q}^{2T} + \frac 15 \begin{cases}
    1 & p+q \equiv 0 \tn{ mod } 4\\
    -1 & p+q \equiv 2 \tn{ mod } 4\\
    0 & \tn{else}
\end{cases}\\
&+\frac 15 \begin{cases}
    1 & p+q \equiv 0 \tn{ mod } 6\\
    -1 & p+q \equiv 4 \tn{ mod } 6\\
    0 & \tn{else}
\end{cases}
+ \frac 15 \begin{cases}
    2 & p+q \equiv 0 \tn{ mod } 10\\
    1 & p+q \equiv 2 \tn{ mod } 10\\
    -1 & p+q \equiv 6 \tn{ mod } 10\\
    -2 & p+q \equiv 8 \tn{ mod } 10\\
    0 & \tn{else}
\end{cases}.
\end{align*}
Every subgroup of $SU(2)$ acts freely on $S^3$. A corollary of \Cref{thm:weyl's law} is that we can discern the order of a discrete group $G \subset SU(2)$ from the spectrum of $\Box_b$ on $G\backslash S^3$. Now, we aim to demonstrate that we can identify the specific group $G$, not just its order. To achieve this, we compare the spectra of all possible quotients of $S^3$ by groups of the same order and perform casework to distinguish them. In particular, we must distinguish the following:
\[2D_{2m},\ \Z/4m\Z\]
\[2T,\ 2D_{12},\ \Z/24\Z\]
\[2O,\ 2D_{24},\ \Z/48\Z\]
\[2I,\ 2D_{60},\ \Z/120\Z.\]
\begin{thm}\label{thm:hear su2}
    Let $G \subset SU(2)$ be a discrete subgroup. We can determine $G$ from the spectrum of $\Box_b$ on $L^2(G \backslash S^3)$.
\end{thm}
\begin{proof}
    We will first compute the multiplicity of $4$, and show that it distinguishes the cyclic groups $\Z/m\Z$ from the non-abelian groups. Then we will compute the multiplicity of $8$, and show that it distinguishes binary dihedral groups from the three exceptional groups.

    First we compute the multiplicity of $4$. The $\mathcal{H}_{p,q}$ spaces with eigenvalue $4$ are $\mathcal{H}_{0,2}$ and $\mathcal{H}_{1,1}$. Let $g \in SU(2)$ with eigenvalues $\mu,\mu^{-1}$. Note that 
    \[\chi_{\mathcal{H}_{0,2}}(g) = \mu^{-2}+1+\mu^2,\tn{ and } \chi_{\mathcal{H}_{1,1}}(g) = \chi_{\mathcal{H}_{0,2}}(g).\]
    So $\dim \mathcal{H}^G_{1,1}=\dim \mathcal{H}^G_{0,2}$ for all $G \subset SU(2)$. For $m > 2$, we have
    \[\dim \mathcal{H}_{1,1}^{\Z/m\Z} = 1.\]
    We compute the dimension 
    \[\dim \mathcal{H}_{1,1}^{2D_{2m}} = \frac 12 \dim \mathcal{H}_{1,1}^{\Z/2m\Z} - \frac 12 = 0.\]
    Then we have the the multiplicity of $4$ in the spectrum of $\Box_b$ on $\Z/4m\Z \backslash S^3$ is $2$, while on $2D_{2m}\backslash S^3$ it is $0$. It follows that this distinguishes between the binary dihedral and cyclic groups.
    
    Note that since $2D_{4}=Q \subset 2T,2O,2I$, it follows that the multiplicity of $4$ in the spectra of $2T \backslash S^3$, $2O \backslash S^3$, $2I \backslash S^3$ is $0$, which distinguishes them from quotients of $S^3$ by cyclic groups.
    
    It remains to distinguish $2T\backslash S^3$ from $2D_{12} \backslash S^3 $, $2O \backslash S^3$ from $2D_{24} \backslash S^3$ and $2I \backslash S^3$ from $2D_{60} \backslash S^3$. We will examine the multiplicity of $8$ to do so. The $\mathcal{H}_{p,q}$ spaces with eigenvalue $8$ are $\mathcal{H}_{3,1},\mathcal{H}_{0,4},\mathcal{H}_{1,2}$. Observe that
    \[\chi_{\mathcal{H}_{1,2}}(g) = \mu^{-3} + \mu^{-1} + \mu + \mu^3,\]
    and
    \[\chi_{\mathcal{H}_{3,1}}(g)=\chi_{\mathcal{H}_{0,4}}(g) = \mu^{-4}+\mu^{-2}+1+\mu^2+\mu^4.\]
    For $k > 3$, we see that
    \[\dim \mathcal{H}_{1,2}^{\Z/k\Z} = 0,\]
    and for $k > 4$,
    \[\dim \mathcal{H}_{3,1}^{\Z/k\Z} = 1.\]
    Also,
    \[\dim \mathcal{H}_{3,1}^{\Z/4\Z} = 3.\]
    From this we see that
    \[\dim \mathcal{H}_{3,1}^{2D_4}=\dim \mathcal{H}_{3,1}^Q = \frac 12 \dim \mathcal{H}_{3,1}^{\Z/4\Z} + \frac 12 = 2 .\]
    and for $k > 2$,
    \[\dim \mathcal{H}_{3,1}^{2D_{2k}}=\frac 12 \dim \mathcal{H}_{3,1}^{\Z/2k\Z} + \frac 12 = 1.\]
    Then
    \[\dim \mathcal{H}_{3,1}^{2T} = \frac 13 \dim \mathcal{H}_{3,1}^Q - \frac 23 = 0.\]
    Which also implies that $\dim \mathcal{H}_{3,1}^{2O} = \dim \mathcal{H}_{3,1}^{2I} = 0$.
    From this, we see that the multiplicity of $8$ in the spectrum is $2$ for $2D_{2k}$ when $k > 2$, while it is $0$ for $2T$, $2O$, $2I$.
\end{proof}

\subsection{Hearing the Fundamental Group}\label{subsec:pi1}
We start by providing a classification of the finite subgroups of $U(2)$ in a format conducive to computing their characters on $\mathcal{H}_{p,q}$ spaces. To facilitate these calculations, we utilize the fact that an element $g \in U(n)$ possesses fixed points if and only if $g$ has an eigenvalue of $1$. 

In this section, we denote an arbitrary cyclic group of order $l$ by $C_l$, while $\Z/m\Z$ exclusively denotes the subgroup of $SU(2)$ generated by the matrix
\[\colvec{e^{2\pi i/m} & 0\\ 0 & e^{-2\pi i/m}}.\]
Similarly, the notations $2D_{2m},Q,2T,2O,2I$ refer to specific subgroups of $SU(2)$ as detailed in \Cref{subsec:su2}.

Our objective is to derive explicit presentations of finite subgroups of $U(2)$ that act freely on $S^3$, as these presentations are essential for character computations.

Consider the mapping $\func{\psi}{SU(2) \times U(1)}{U(2)}$ defined by multiplication, where $U(1)$ is identified with the center $Z(U(2))$. This mapping is a double cover of $U(2)$. 
Note that the natural action of $SU(2) \times U(1)$ on $S^3$ factors through $\psi$. A finite subgroup $G \subset U(2)$ lifts to a group $\psi^{-1}(G)$ of twice the order, containing 
\[N = \ker \psi = \left\{\left(\colvec{-1 & 0\\ 0 & -1},-1 \right),\left( \colvec{1 & 0\\ 0 & 1}, 1\right)\right\}.\]

Additional clarification of \Cref{prop:subgrps u(2)} will be given after to not clutter up the statement.
\begin{prop}\label{prop:subgrps u(2)}
    There are six families of non-abelian subgroups of $U(2)$ acting freely on $S^3$ up to conjugation, these are
    \begin{enumerate}
        \item $2I \times C_l$, $l$ coprime to $30$. 
        \item $2O \times C_l$, $l$ coprime to $6$.
        \item $2T \times C_l$, $l$ coprime to $6$.
        \item $2D_{2m} \times C_l$, $l$ coprime to $2m$.
        \item $Q \rtimes C_{18l}/N$, $l$ odd.
        \item $\Z/2m\Z \rtimes C_{4l}/N$, $m$ odd, $l$ even, $m$ coprime to $l$.
    \end{enumerate}
    The only abelian subgroups of $U(2)$ acting freely on $S^3$ are cyclic groups $C_l$. There are $1+ \frac{\phi(l)}{2}$ of these up to conjugation for $l>2$, where $\phi$ is the Euler totient function. For $l=2$, there is only one.
\end{prop}
Types 1-4 of the form $G \times C_l$ are generated by $G \subset SU(2)$ and additionally
\[\colvec{e^{2\pi i/l} & 0\\ 0 & e^{2\pi i/l}}.\]
Type 5 is generated by $Q$ along with
\[\frac{e^{\pi i/(9l)}}2\colvec{1+i & -1+i\\ 1+i & 1-i}.\]
Type 6 is generated by $\Z/2m\Z$ and
\[e^{\pi i /(2l)}\colvec{0 & -1\\ 1 & 0}.\]
This classification is known but we could not locate a source for the given presentations in the literature.
\begin{proof}
We will examine the finite subgroups of $SU(2)\times U(1)$ of the form $\psi^{-1}(G)$, where $G$ acts freely on $S^3$. Let $p_1,p_2$ be the projections onto the first and second factors respectively. Let $G$ be a finite subgroup of $SU(2)\times U(1)$. Then $p_1(G)$ and $G \cap \ker p_2 \cong p_1(G \cap \ker p_2)$ are both finite subgroups of $SU(2)$. Then
\[\frac{p_1(G)}{p_1(G \cap \ker p_2)} \cong \frac{G/(G \cap \ker p_1)}{((G \cap \ker p_2)(G \cap \ker p_1))/(G\cap \ker p_1)} \cong\] \[\cong \frac{G}{(G \cap \ker p_2)(G \cap \ker p_1)} \cong \frac{p_2(G)}{p_2(G \cap \ker p_1)}.\]
So the quotient of $p_1(G)$ by $p_1(G \cap \ker p_2)$ is cyclic. We will classify $G$ based on $p_1(G)$.

\subsubsection*{{\bf Case I: $p_1(G) \cong 2I$}}\

The normal subgroups of $2I$ are $1,\{\pm 1\}, 2I$. Only $2I/2I$ is cyclic. In this case, $G$ is a member a family of subgroups of $SU(2) \times U(1)$ of the form $2I \times C_k$. Note there is only one up to conjugacy for each $k$. When $k=2l$ is even, the group $2I \times C_{2l}$ contains the kernel of $\psi$. When $l$ is odd, $\psi(2I \times C_{2l}) \cong 2I \times C_l$. This group is generated by $2I$ with 
\[\colvec{e^{2\pi i/l} & 0 \\ 0 & e^{2\pi i/l}}.\]
We observe that $\psi(2I \times C_{2l})$ acts freely on $S^3$ if and only if $l$ is relatively prime to $2*3*5=30$. We will refer to this family of groups and those conjugate to them by $2I \times C_l$.

\subsubsection*{{\bf Case II: $p_1(G) \cong 2O$}}\ 

The normal subgroups of $2O$ are $1,\{\pm 1\}, Q, 2T, 2O$, where $Q$ is the quaternion group. Only $2O/2O$ and $2O/2T$ are cyclic. If $G \cap \ker p_2 \cong 2O$, then $G$ is of the form $2O \times C_k$. This contains the kernel of $\psi$ when $k=2l$ is even. When $l$ is odd, $\psi(2O \times C_{2l}) \cong 2O \times C_l$. This group is generated by $2O$ with 
\[\colvec{e^{2\pi i/l} & 0 \\ 0 & e^{2\pi i/l}}.\]
We observe that $\psi(2O \times C_{2l})$ acts freely on $S^3$ if and only if $l$ is relatively prime to $2 \cdot 3=6$. We will refer to these by $2O \times C_l$.

If $G \cap \ker p_2 \cong 2T$, then $G \cap \ker p_1 \cong C_k$, while $p_2(G) \cong C_{2k}$. This contains the kernel of $\psi$ when $k=2l$ is even. Consider the homomorphism $C_{4l} \to G$, given by 
\[y \mapsto \left(g, e^{\pi i / (2l)}\right).\]
Where $g$ is an order $4$ element of $2O \smallsetminus 2T$, and $y$ is a generator of $C_{4l}$. Therefore, the short exact sequence
\[
\begin{tikzcd}
    0 \arrow[r] & 2T \arrow[r] & G \arrow[r] & C_{4l} \arrow[r] & 0
\end{tikzcd}.
\]
splits, so $G$ is a semidirect product $2T \rtimes C_{4l}$. When $l$ is odd, $(g,i) \in G$ for $g$ being an order $4$ element of $2O \smallsetminus 2T$. So $g$ has eigenvalues $\pm i$. Then $ig \in U(2)$ has eigenvalues $\pm 1$. When $l$ is even, we have $(g',i) \in G$ for $g'$ being an order $4$ element of $2T$. So $\psi(G)$ does not act freely on $S^3$.

\subsubsection*{{\bf Case III: $p_1(G) \cong 2T$}}\ 

The normal subgroups of $2T$ are $1, \{\pm 1\}, Q, 2T$. Only $2T/2T$ and $2T/Q$ are cyclic. If $G \cap \ker p_2 \cong 2T$, then $\psi(G)$ is a member of a family of groups $2T \times C_l$, where $l$ is relatively prime to $6$. We will refer to these by $2T \times C_l$.

If $G \cap \ker p_2 \cong Q$, then there are two superficially different possibilities for $G$. We have $2T = Q \sqcup aQ \sqcup bQ$. Each of the cosets $aQ$, $bQ$ contains $4$ elements of order $3$ and $4$ elements of order $6$. The elements of order $3$ are all conjugate in $SU(2)$. So $(g,e^{2\pi i/(3k)}) \in G$ where $g \in 2T \smallsetminus Q$ is an order $3$ element. This shows $G$ is a semidirect product $Q \rtimes C_{3k}$. This contains the kernel of $\psi$ when $k$ is even. If $3$ does not divide $k$, then we see
$g^ke^{2\pi ik/(3k)}$ has eigenvalues $1$ and $e^{4\pi i/3}$. Then $\psi(G)$ does not act freely on $S^3$. Let $6l=k$. We see that $Q \rtimes C_{18l}$ contains the index $3$ subgroup $Q \times C_{6l}$, and the elements outside of this subgroup do not have fixed points in their action on $S^3$. Note that if $l$ is odd, then $\psi(Q \rtimes C_{18l})$ acts freely on $S^3$. We will refer to this family of groups by $Q \rtimes C_{18l}/N$, where $N = \ker \psi$.

\subsubsection*{{\bf Case IV: $p_1(G) \cong 2D_{2m}$}}\ 

A convienient presentation for $2D_{2m}$ is $2D_{2m}=\<a,x\mid a^mx^2=x^4=ax^{-1}ax=1\>$. The normal subgroups of $2D_{2m}$ are $1, 2D_{2m}, \<a^k\>$ for $k | 2m$. When $m$ is even, there are two more: $\<a^2,x\>,\ \<a^2,ax\>$. When $m$ is odd, only $2D_{2m}/2D_{2m},\ 2D_{2m}/\<a\>$ and $ 2D_{2m}/\<a^2\>$ are cyclic. When $m$ is even, only $2D_{2m}/2D_{2m}$ and $2D_{2m}/\<a\>$ are cyclic. The case $G \cap \ker p_2 = 2D_{2m}$ is very similar to previous cases. We will refer to this family of groups by $2D_{2m} \times C_l$.

If $G \cap \ker p_2 = \<a\>$, then $G \cap \ker p_1 \cong C_k$, while $p_2(G) \cong C_{2k}$. This contains the kernel when $k=2l$ is even. Consider the homomorphism $C_{4l} \rightarrow G$ given by
\[y \mapsto (x, e^{\pi i/(2l)}).\]
Then $G$ is a seimdirect product $\Z/2m\Z \rtimes C_{4l}$. Note $(x^l,i) \in G$. If $l$ is odd, then $ix^l$ has eigenvalues $\pm 1$. So $\psi(G)$ does not act freely on $S^3$. If $l$ and $m$ are even, then $(a^{m/2},i) \in G$. Since $ia^{m/2}$ has eigenvalues $\pm 1$, $\psi(G)$ does not act freely on $S^3$. Suppose $m$ is odd and $l$ is even. First consider elements of the form
\[(a^jx, e^{\pi i(2o+1)/(2l)})\]
for $o=0,\ldots , 2l-1$. We see $e^{\pi i(2o+1)/(2l)}a^jx$ does not have any eigenvalues equal to $1$. The remaining elements in $G$ are a subgroup $\Z/2m\Z \times C_{2l}$. We observe this is generated by 
\[(a,1),(1,e^{\pi i/l}).\]
Suppose $c=\tn{gcd}(l,m)>1$. Then we see that $a^{m/c}(e^{\pi i/l})^{l/c}$ has eigenvalues $1, e^{2\pi i/c}$. Then $\psi(G)$ does not act freely on $S^3$. If $m,l$ are coprime, then note $\Z/2m\Z \times C_{2l} \cong C_2 \times C_{2lm}$. The element
\[(a^2,e^{\pi i/l})\]
has order $2lm$.
We observe that $\psi(\Z/2m\Z\times C_{2l}) \cong C_{2lm}$,
and is generated by a matrix conjugate to
\[\colvec{e^{\frac{\pi i(m+2l)}{ml}} & 0\\ 0 & e^{\frac{\pi i(m-2l)}{ml}}}.\]
This evidently acts freely on $S^3$. So $\psi(\Z/2m\Z \rtimes C_{4l})$ acts freely on $S^3$ when $m$ is odd, $l$ is even, and $l$ is coprime to $m$. We refer to this family of groups by $\Z/2m\Z \rtimes C_{4l}/N$.

When $m$ is odd and $G \cap \ker p_2 = \<a^2\>$, we have that $G \cap \ker p_1 \cong C_k$, while $p_2(G) \cong C_{4k}$. There are two possibilities, $(x,e^{\pi i /(2k)}) \in G$ or $(x^3,e^{\pi i/(2k)}) \in G$. Since $x,x^3$ are conjugate in $SU(2)$, the two possibilities are conjugate. Then we have the homomorphism $C_{4k} \rightarrow G$, 
\[y \mapsto (x,e^{\pi i /(2k)})\]
Then $G$ is a semidirect product $\Z/m\Z \rtimes C_{4k}$. $G$ contains the kernel of $\psi$ if and only if $k$ is odd. Then $(x^k, e^{\pi i k/(2k)}) \in G$. But this has eigenvalues $\pm 1$. So $G$ does not act freely on $S^3$.

\subsubsection*{{\bf Case V: $p_1(G) \cong \Z/m\Z$}}\ 

The abelian case is far simpler.
\begin{thm}[The $p^2$-condition \cite{smithp^2}]
    Let $p$ be a prime. The group $G = (\Z/p\Z)^2$ does not have a continuous free action on any sphere.
\end{thm}
The $p^2$-condition shows that the only abelian groups that act freely on $S^3$ are cyclic groups. Unlike the non-abelian cases, there generally are more than one of these up to conjugation. The exact number is simple to compute \cite{GLR}.
\end{proof}

To reiterate, we have seven types of groups:
\begin{enumerate}
        \item $2I \times C_l$, $l$ coprime to $30$.
        \item $2O \times C_l$, $l$ coprime to $6$.
        \item $2T \times C_l$, $l$ coprime to $6$.
        \item $2D_{2m} \times C_l$, $l$ coprime to $2m$.
        \item $Q \rtimes C_{18l}/N$, $l$ odd.
        \item $\Z/2m\Z \rtimes C_{4l}/N$, $m$ odd, $l$ even, $m$ coprime to $l$.
        \item Cyclic groups $C_l$, of which there are $1+ \frac{\phi(l)}{2}$ of these up to conjugation for $l>2$. For $l=2$, there is only one.
\end{enumerate}

We note that type 3 and 5 are both central extensions of the tetrahedral group $A_4$, while types 4 and 6 are both central extensions of dihedral groups.

    We will show that one can distinguish between the seven infinite families listed above, and that $l,m$ can be determined for type 4 and 6. For the other types, the order can be determined by \Cref{thm:weyl's law}, and this gives which member of the family it is.

The proof below occupies the remainder of this section.
\begin{proof}(Of \Cref{thm:hear pi1})\\
    This proof proceeds by cases. We compute the multiplicities of $4$, $8$, $12$ and $24$ in the spectrum of $\Box_b$ on certain sphere quotients. The comparisons are as follows:\\
    \begin{itemize}
        \item The multiplicity of $4$ is positive for $G$ cyclic, and it is zero for $G$ non-abelian.
        \item No group of the form $Q \rtimes C_{18l}/N$ has the same order as a groups of the form $2T\times C_{l'}$, $2O \times C_{l'}$ or $2I \times C_{l'}$. By \Cref{thm:weyl's law}, they are not isospectral.
        \item No group of the form $2O \times C_l$ has the same order as groups of the form $2T \times C_{l'}$ or $2I \times C_l$. By \Cref{thm:weyl's law}, they are not isospectral.
        \item The multiplicity of $24$ differs for groups of the form $2I \times C_l$ and $2T \times C_{5l}$.
        \item The multiplicity of $12$ for $2D_{2m} \times C_l$ differs from the multiplicity of $12$ for $2T \times C_{l'}$, $2O \times C_{l'}$, $2I \times C_{l'}$ and $Q \rtimes C_{18l'}/N$.
        \item The multiplicity of $12$ for $C_{2m} \rtimes C_{4l}/N$ differs from the multiplicity of $12$ for $2T \times C_{l'}$, $2O \times C_{l'}$, $2I \times C_{l'}$ and $Q \rtimes C_{18l'}/N$, except for $2T$ compared to $\Z/6\Z \rtimes C_8/N$. The multiplicity of $24$ differs for $2T$ and $\Z/6\Z \rtimes C_8/N$.
        \item The multiplicity of $8$ is positive for $2D_{2m}$, and zero for $2D_{2m'}\times C_l$, $l > 1$ and $C_{2m''} \rtimes C_{4l'}/N$. To compare $2D_{2m} \times C_l$, $l > 1$, with $C_{2m'} \rtimes C_{4l'}/N$, we show the multiplicity of $4r$ differs, for $r$ some sufficiently large prime such that $r \cong 1 \tn{ mod } 4ll'$. The same computation works to compare $2D_{2m} \times C_l$ with $2D_{2m'} \times C_{l'}$, for $1 < l < l'$. This also works for comparing $C_{2m} \rtimes C_{4l}/N$ with $C_{2m'} \times C_{4l'}/N$ for $l < l'$.
    \end{itemize}
    Recall $\func{\psi}{SU(2)\times U(1)}{U(2)}$ is the double cover given by multiplication with $U(1)$ identified with the center of $U(2)$, and $N = \ker \psi$. Factoring through $\psi$, subgroups of $SU(2)\times U(1)$ act on $S^3$. We observe that
    \[\dim \mathcal{H}_{p,q}^{\psi^{-1}(G)} = \dim \mathcal{H}_{p,q}^{G}\]
    Where $G$ is a finite subgroup of $U(2)$.
    
    We begin by computing the characters of the non-abelian groups. Let $\lambda, \zeta$ be unit norm complex numbers. We compute
    \[\chi_{H_{p,q}}\left(\colvec{\lambda & 0\\0 & \overline{\lambda}}\colvec{\zeta & 0\\ 0 & \zeta} \right) = (\zeta\overline{\lambda})^q(\zeta\lambda)^{-p}\sum_{j=0}^{p+q}\left( \frac{\zeta \lambda}{\zeta \overline{\lambda}} \right)^j=\zeta^{q-p}\chi_{H_{p,q}}\left( \colvec{\lambda & 0\\0 & \overline{\lambda}} \right).\]

    For groups of the form $G \times C_{2l} \subset SU(2) \times U(1)$, where $G \subset SU(2)$ is a finite subgroup, and $C_{2l} \subset U(1)$, we compute
    \[\dim \mathcal{H}_{p,q}^{\psi(G \times C_{2l})} = \dim \mathcal{H}_{p,q}^{G \times C_{2l}} = \frac 1{2l} \sum_{j=0}^{2l-1}(e^{\pi i j/l})^{j(q-p)}\dim \mathcal{H}_{p,q}^{G} = \begin{cases}
        \dim \mathcal{H}_{p,q}^{G} & q \equiv p \tn{ mod } 2l\\
        0 & \tn{else}
    \end{cases}.\]
    Recall that $\dim \mathcal{H}_{p,q}^{\Z/2m\Z} = 0$ when $p+q$ is odd, which implies that $\dim\mathcal{H}_{p,q}^{2D_{2m}}=\dim\mathcal{H}_{p,q}^{2T}=\dim\mathcal{H}_{p,q}^{2O}=\dim\mathcal{H}_{p,q}^{2I}=0$ when $p+q$ is odd. For $G$ one of $2I,2O,2T,2D_{2m}$, we have
     \[\dim \mathcal{H}_{p,q}^{G \times C_{l}} = \begin{cases}
        \dim \mathcal{H}_{p,q}^{G} & q \equiv p \tn{ mod } l\\
        0 & \tn{else}
    \end{cases},\]
    since $l$ is odd. Recall that $Q \rtimes C_{18l}$ is generated by $Q$ and $(g,e^{\pi i/(9l)}) \in SU(2) \times U(1)$ with $g \in 2T$ an order $3$ element. Then we have that $gQ$ has four elements with order $3$ and four elements with order $6$. The remaining coset $g^2Q$ has four elements with order $3$ and four elements with order $6$. Let $\zeta = e^{\pi i/(9l)}$. Then we compute:
    {\begin{align*}
        \dim \mathcal{H}_{p,q}^{Q \rtimes C_{18l}} &= \frac 1{18l} \sum_{j=0}^{6l-1}(\zeta^{q-p})^{3j}\dim\mathcal{H}_{p,q}^Q +\\ &+ \frac 1{36l}\sum_{k=1}^2\sum_{j=0}^{6l-1}\zeta^{k(q-p)}(\zeta^{q-p})^{3j}\left[ e^{\pi i (p+q)/3} + e^{-2\pi i(p+q)/3} \right]\sum_{m=0}^{p+q}e^{-2\pi im/3}\\
        &= \frac 13\begin{cases}
            \dim H_{p,q}^Q & p\equiv q \tn{ mod } 6l\\
            0 & \tn{else}
        \end{cases}\ \  + \\
        &+ \frac 16 \left[\begin{cases}
            -1 & q-p \cong 6l \tn{ mod } 18l\\
            -1 & q-p \cong 12l \tn{ mod } 18l\\
            2 & q-p \cong 0 \tn{ mod } 18l\\
            0 & \tn{else}
        \end{cases} \right]\begin{cases}
            2 & p+q\equiv 0 \tn{ mod } 6\\
            -2 & p+q\equiv 4 \tn{ mod } 6\\
            0 & \tn{else}
        \end{cases}.
    \end{align*}}
    Recall that $\Z/2m\Z \rtimes C_{4l}$ is generated by $\Z/2m\Z \subset SU(2)$ and $(x,e^{\pi i/2l})$, where $x \in 2D_{2m}$ is an order $4$ element. Every element of $x\Z/2m\Z \subset 2D_{2m}$ has order $4$, so they all have eigenvalues $\pm i$. Let $\zeta = e^{\pi i/2l}$. We compute:
    {\begin{align*}
        \dim \mathcal{H}_{p,q}^{\Z/2m\Z \rtimes C_{4l}} &= \frac 1{4l}\sum_{j=0}^{2l-1}(\zeta^{q-p})^{2j}\dim \mathcal{H}_{p,q}^{\Z/2m\Z} + \frac 1{4l}\sum_{j=0}^{2l-1}\zeta^{q-p}(\zeta^{q-p})^{2j}i^{-p-q}\sum_{k=0}^{p+q}(-1)^k\\
        &= \frac 12 \begin{cases}
            \dim \mathcal{H}_{p,q}^{\Z/2m\Z} & p \equiv q \tn{ mod } 2l\\
            0 & \tn{else}
        \end{cases}\ \ + \\
        &+ \frac 12\left[\begin{cases}
            -1 & q-p \equiv 2l \tn{ mod } 4l\\
            1 & q-p \equiv 0 \tn{ mod } 4l\\
            0 & \tn{else}
        \end{cases}\right]\begin{cases}
            1 & p+q \equiv 0 \tn{ mod } 4\\
            -1 & p+q \equiv 2 \tn{ mod } 4\\
            0 & \tn{else}
        \end{cases}.
    \end{align*}}
    We first compute the multiplicity of $4$ for $\Z/2m\Z \rtimes C_{4l}/N$. The $\mathcal{H}_{p,q}$ spaces with this eigenvalue are $\mathcal{H}_{1,1}$ and $\mathcal{H}_{0,2}$. Since $l$ is even, we observe that
    \[\dim \mathcal{H}_{0,2}^{\Z/2m\Z \rtimes C_{4l}} = 0.\]
    Since $\dim \mathcal{H}_{1,1}^{\Z/2m\Z} = 1$ for $m > 2$, we have
    \[\dim \mathcal{H}_{1,1}^{\Z/2m\Z \rtimes C_{4l}} = \frac 12 (1) + \frac 12 (1)(-1) = 0.\]
    We also have that the multiplicity of $4$ in the spectrum of $\Box_b$ is $0$ for quotients by non-abelian subgroups of $SU(2)$. Each of the non-abelian groups of types 1-5 in \Cref{prop:subgrps u(2)} contains a non-abelian subgroup of $SU(2)$. So we see that the multiplicity of $4$ in the spectrum of $\Box_b$ is $0$ for all sphere quotients by a non-abelian finite subgroup of $U(2)$. Since $\frac{\partial}{\partial z_2}\frac{\partial}{\partial \overline{z}_2}|z|^{-2} \in \mathcal{H}_{1,1}$ is invariant under diagonal unitary matrices, we have that the multiplicity of $4$ for any lens space is non-zero. So this distinguishes the abelian quotients from the non-abelian quotients.\\
    
    We can distinguish between groups of different order by \Cref{thm:weyl's law}. In particular, we note that a group of the form $2I \times C_l$ does not have the same order as a group of the form $2O \times C_{l'}$. Otherwise, $5l = 2l'$, a contradiction as $l$ must be odd. Similarly, a group of the form $2O \times C_l$ does not have the same order as a group of the form $2T \times C_{l'}$, as otherwise $l' = 2l$. A group of the form $Q \rtimes C_{18l}/N$ does not have the same order as a groups of the form $2T \times C_{l'}$, $2O \times C_{l'}$, $2I \times C_{l'}$, as any of these would imply $3 \mid l'$.

    To distinguish between $2I \times C_l$ and $2T \times C_{5l}$, we analyze the multiplicity of $24$. The $\mathcal{H}_{p,q}$ spaces with eigenvalue $24$ for $\Box_b$ are $\mathcal{H}_{11,1},\mathcal{H}_{5,2},\mathcal{H}_{3,3},\mathcal{H}_{2,4},\mathcal{H}_{1,6},\mathcal{H}_{0,12}$. Since
    \[\dim \mathcal{H}_{5,2}^{\Z/4\Z} = \dim \mathcal{H}_{1,6}^{\Z/4\Z} = 0,\]
    the dimensions of these spaces for $2I \times C_l$ and $2T \times C_{5l}$ are zero.
    We compute
    \begin{alignat*}{2}
        \dim \mathcal{H}_{3,3}^{\Z/4\Z} = \dim\mathcal{H}_{2,4}^{\Z/4\Z} &= \frac 14\sum_{j=0}^3(i^{-6})^j\sum_{k=0}^{6}(i^{2k})^j &=3\\
        \dim \mathcal{H}_{3,3}^{Q} &= \frac 12 (3) + \frac 12 (-1) &= 1\\
        \dim \mathcal{H}_{3,3}^{2T} &= \frac 13 (1) + \frac 13 (2) &= 1\\
        \dim \mathcal{H}_{3,3}^{2I} &= \frac 15 (1) + \frac 15 (-1) + \frac 15 (1) + \frac 15 (-1) &= 0\\
        \dim \mathcal{H}_{0,12}^{\Z/4\Z} &= \frac 14\sum_{j=0}^3(i^{-12})^j\sum_{k=0}^{12}(i^{2k})^j &= 7\\
        \dim \mathcal{H}_{0,12}^{Q} &= \frac 12 (7) + \frac 12(1) &= 4\\
        \dim \mathcal{H}_{0,12}^{2T} &= \frac 13 (4) + \frac 13 (2) &= 2\\
        \dim \mathcal{H}_{0,12}^{2I} &= \frac 15 (2) + \frac 15 (1) + \frac 15 (1) + \frac 15 (1) &= 1.
    \end{alignat*}
    Recall that $l$ is relatively prime to $30$. Then we see if $l > 1$, $12-0=12$ and $1-11=-10$ are not congruent to $0$ modulo $2l$. So the multiplicity of $24$ for $2I \times C_l$ is zero when $l > 1$. For $2T \times C_{5l}$, it is at least $1$. When $l=1$, we see that the multiplicity of $24$ for $2I$ is $2$. We see that
    \[\dim \mathcal{H}_{0,12}^{2T\times C_{5l}} = \dim \mathcal{H}_{2,4}^{2T\times C_{5l}} = 0,\]
    but
    \[\dim \mathcal{H}_{3,3}^{2T\times C_{5l}} = 1,\ \dim \mathcal{H}_{11,1}^{2T\times C_{5l}} = 2.\]
    So the multiplicity of $24$ for $2T \times C_{5l}$ is $3$.

    Next we compute the multiplicity of $12$, to distinguish groups of the types $2D_{2m} \times C_l$ and $\Z/2m\Z \rtimes C_{4l}/N$ from the other types. The $\mathcal{H}_{p,q}$ spaces that have eigenvalue $12$ are $\mathcal{H}_{5,1},\mathcal{H}_{2,2},\mathcal{H}_{1,3},\mathcal{H}_{0,6}$. From previous computations, we have
    \[\dim \mathcal{H}_{2,2}^{\Z/4\Z} = 3,\ \dim \mathcal{H}^Q_{2,2}=2,\ \dim \mathcal{H}_{2,2}^{2T} = \dim \mathcal{H}_{2,2}^{2O}=\dim \mathcal{H}_{2,2}^{2I} = 0.\]
    and for $m > 2$,
    \[\dim \mathcal{H}_{2,2}^{\Z/2m\Z} = \dim \mathcal{H}_{2,2}^{2D_{2m}} = 1.\]
    Then we have
    \[\dim \mathcal{H}_{2,2}^{\Z/2m\Z \rtimes C_{4l}} = \frac 12 (1) + \frac 12 (1)(1) = 1,\]
    and
    \[\dim \mathcal{H}_{1,3}^{\Z/2m\Z \rtimes C_{4l}} = 0.\]
    We compute
    \[\dim \mathcal{H}_{0,6}^{\Z/6\Z} = 3,\ \dim \mathcal{H}_{0,6}^{2D_6} = \frac 12 (3) + \frac 12(-1) = 1,\]
    and for $m > 3$,
    \[\dim \mathcal{H}_{0,6}^{\Z/2m\Z} = 1,\ \dim \mathcal{H}_{0,6}^{2D_{2m}} = \frac 12 (1) + \frac 12(-1) = 0.\]
    We compute
    \begin{alignat*}{2}
        \dim \mathcal{H}_{0,6}^{2O} &= \frac 12 (1) + \frac 12 (-1) &= 0\ \\
        \dim \mathcal{H}_{2,2}^{Q \rtimes C_{18l}} &= \frac 13 (2) + \frac 16 (2) (-2) &= 0\ \\
        \dim \mathcal{H}_{0,6}^{Q \rtimes C_{18}} &= \frac 13 (1) + \frac 16 (-1) (2) &= 0.
    \end{alignat*}
    For $l > 1$, $\dim \mathcal{H}_{0,6}^{Q \rtimes C_{18l}}=0$, and $\dim \mathcal{H}_{5,1}^{Q \rtimes C_{18l}}=\dim \mathcal{H}_{1,3}^{Q \rtimes C_{18l}}=0$ for all $l$. So we see that the multiplicity of $12$ is at least $1$ for $2D_{2m} \times C_l$ and $\Z/2m\Z \rtimes C_{4l}/N$, while it is zero for $2I \times C_l, 2O \times C_l, Q \rtimes C_{18l}/N$. We see that for $l>1$, 
    \[\dim \mathcal{H}_{5,1}^{2T \times C_l} = \dim \mathcal{H}_{0,6}^{2T \times C_l} = 0.\]
    Since $l$ is relatively prime to $6$. For $l=1$, we have
    \[\dim \mathcal{H}_{5,1}^{2T} = \dim \mathcal{H}_{0,6}^{2T} = 1.\]
    So the multiplicity is $2$ for $2T$. The orders of $2D_{2m}\times C_l$ and $\Z/2m\Z \rtimes C_{4l}/N$ are both $4ml$. From \Cref{thm:weyl's law}, we know the spectra differ if $4ml \neq 24$. The groups that have this order are $2D_{12},Q \times C_3, \Z/6\Z \rtimes C_8/N$. From \Cref{thm:hear su2}, the spectra of $2D_{12}$ and $2T$ differ. We see that
    \[\dim\mathcal{H}_{2,2}^{Q \times C_3} = 2,\ \dim\mathcal{H}_{1,3}^{Q \times C_3} = \dim\mathcal{H}_{5,1}^{Q \times C_3}=0,\ \dim\mathcal{H}_{0,6}^{Q \times C_3}=1\]
    So the multiplicity of $12$ for $Q \times C_3$ is $3$. For $\Z/6\Z \rtimes C_8/N$, we need to compute the multiplicity of $24$. For $2T$, it is $6$. We compute
    \begin{alignat*}{2}
        \dim \mathcal{H}_{3,3}^{\Z/6\Z \rtimes C_{8}} &= \frac 12 (3) + \frac 12 (1)(-1) &= 1\ \\
        \dim \mathcal{H}_{0,12}^{\Z/6\Z} & &= 5\ \\
        \dim \mathcal{H}_{0,12}^{\Z/6\Z \rtimes C_{8}} &= \frac 12 (5) + \frac 12 (-1)(1) &= 2.
    \end{alignat*}
    Since $\dim \mathcal{H}_{11,1}^{\Z/6\Z \rtimes C_{8}}=\mathcal{H}_{2,4}^{\Z/6\Z \rtimes C_{8}} = 0$, we see that the multiplicity of $24$ is $3$ for $\Z/6\Z \rtimes C_{8}/N$.

    Now we distinguish $2D_{2m} \times C_l$ from $\Z/2m'\Z \rtimes C_{4l'}/N$, with $ml=m'l'$. Note that $l$ is odd and $l'$ is even, so $l \neq l'$, and $m\neq m'$. First suppose $l > 1$. Note there are infinitely many primes $r$ such that $r \equiv 1 \tn{ mod } 4ll'$. Pick one such that
    \[r > \frac{2}{\left| \frac{1}{m}- \frac 1{m'}\right|}.\]
    The $\mathcal{H}_{p,q}$ spaces with eigenvalue $4r$ are $\mathcal{H}_{2r-1,1},\mathcal{H}_{r-1,2},\mathcal{H}_{1,r},\mathcal{H}_{0,2r}$. We observe that 
    \begin{align*}
        1-(2r-1) &\equiv 0 \tn{ mod } 4ll'\\
        2-(r-1) &\equiv 2 \tn{ mod } 4ll'\\
        r-1&\equiv 0 \tn{ mod } 4ll'\\
        2r-0 &\equiv 2 \tn{ mod } 4ll'
    \end{align*}
    Since $l>1$ and $l$ is odd, $2 \not \equiv 0 \tn{ mod } l$. Since $2l' \geq 4$, then $2 \not \equiv 0 \tn{ mod } 2l'$. So then we see
    \[\dim \mathcal{H}_{r-1,2}^{2D_{2m} \times C_l} = \dim \mathcal{H}_{r-1,2}^{\Z/2m'\Z \rtimes C_{4l'}/N}=\dim \mathcal{H}_{0,2r}^{2D_{2m} \times C_l} = \dim \mathcal{H}_{0,2r}^{\Z/2m'\Z \rtimes C_{4l'}/N}=0.\]
    We have
    \begin{align*}
        \dim\mathcal{H}_{2r-1,1}^{\Z/2m\Z} &= 2\left\lfloor \frac{2r-1}{2m} \right\rfloor + 1\\
        \dim\mathcal{H}_{1,r}^{\Z/2m\Z} &=  2\left\lfloor \frac{r}{2m} \right\rfloor + 1.
    \end{align*}
    Then we compute
    \begin{alignat*}{2}
        \dim\mathcal{H}_{2r-1,1}^{2D_{2m}\times C_l} &= \frac 12 \left(2\left\lfloor \frac{2r-1}{2m} \right\rfloor + 1\right) - \frac 12 &= \left\lfloor \frac{2r-1}{2m} \right\rfloor\ \\
        \dim\mathcal{H}_{1,r}^{2D_{2m}\times C_l} &= \frac 12 \left(2\left\lfloor \frac{r}{2m} \right\rfloor + 1\right) - \frac 12 &= \left\lfloor \frac{r}{2m} \right\rfloor\ \\
        \dim\mathcal{H}_{2r-1,1}^{\Z/2m'\Z \rtimes C_{4l'}} &= \frac 12 \left( 2\left\lfloor \frac{2r-1}{2m'} \right\rfloor + 1 \right) + \frac 12 (1) (-1) &= \left\lfloor \frac{2r-1}{2m'} \right\rfloor\ \\
        \dim\mathcal{H}_{1,r}^{\Z/2m'\Z \rtimes C_{4l'}} &= \frac 12 \left( 2\left\lfloor \frac{r}{2m'} \right\rfloor + 1 \right) + \frac 12 (1) (-1) &= \left\lfloor \frac{r}{2m'} \right\rfloor.
    \end{alignat*}
    Suppose $m > m'$. Then we compute
    \begin{align*}
        \left\lfloor \frac{r}{2m'} \right\rfloor + \left\lfloor \frac{2r-1}{2m'} \right\rfloor-\left\lfloor \frac{r}{2m} \right\rfloor-\left\lfloor \frac{2r-1}{2m} \right\rfloor &\geq \frac{r}{2m'}+\frac{2r}{2m'}-\frac{r}{2m}-\frac{2r}{2m}-\left(\frac 1{m'} - \frac 1{m} \right) -2\\
        &> \frac 32 r\left(\frac{1}{m'}-\frac 1{m} \right) - 3\\
        &> 0.
    \end{align*}
    If $m' > m$, the above computation works with $m$ and $m'$ swapped. So we see that the multiplicity of $4r$ differs for $2D_{2m} \times C_l$ and $\Z/2m'\Z \rtimes C_{4l'}/N$.
    Suppose $l=1$. Then we will examine the multiplicity of $8$. The $\mathcal{H}_{p,q}$ spaces with eigenvalue $8$ are $\mathcal{H}_{3,1},\mathcal{H}_{1,2},\mathcal{H}_{0,4}$. We see from \Cref{thm:hear su2} that the multiplicity of $8$ for $2D_{2m}$ is at least $2$. Since $2l' \geq 4$, we see that 
    \[\dim \mathcal{H}_{3,1}^{\Z/2m'\Z \rtimes C_{4l'}} = \dim \mathcal{H}_{1,2}^{\Z/2m'\Z \rtimes C_{4l'}}=0.\]
    When $l' > 2$, then we have
    \[\dim \mathcal{H}_{0,4}^{\Z/2m'\Z  \rtimes C_{4l'}}=0.\]
    Otherwise,
    \[\dim \mathcal{H}_{0,4}^{\Z/2m'\Z  \rtimes C_{8}}=\frac 12 (1) + \frac 12 (-1) (1) = 0.\]
    So the multiplicity of $8$ is $0$.\\

    A nearly identical proof to the one above works for comparing $\Z/2m\Z \rtimes C_{4l}/N$ with $\Z/2m'\Z \rtimes C_{4l'}/N$ when $ml=m'l'$ and $l\neq l'$. This also works for comparing $2D_{2m} \times C_{l}$ with $2D_{2m'}\times C_{l'}$ when $ml=m'l'$ and $1 < l < l'$ (recall that $l,l'$ are odd). When $l=1 < l'$, we observe that $l' \neq 2,4$, so then we have
    \[\dim \mathcal{H}_{3,1}^{2D_{2m'}\times C_{l'}}=\dim \mathcal{H}_{0,4}^{2D_{2m'}\times C_{l'}}=0.\]
    So the multiplicity of $8$ for $2D_{2m'}\times C_{l'}$ is zero, while the multiplicity of $8$ for $2D_{2m}$ is at least $2$.
\end{proof}

\section{Sobolev Estimates for the Complex Green's Operator}\label{sec:Sobolev}
In this section, we prove Sobolev estimates on sphere quotients. We use the following characterization of the Sobolev space as the space of Bessel potentials (see \cite{stein-singular-integrals} for a reference) to demonstrate that the complex Green's operator $\mathcal{G}$ gives one weak derivative.
\[H^s(M) = \{f \in L^2(M) : (I + \Delta_{M})^{s/2}f \in L^2(M)\}\]
For more background on Sobolev estimates, we recommend the reader consult \cite{Chen-Shaw}.
The eigenspaces of $\Delta_{S^{2n-1}}$ are the spaces of homogenous harmonic polynomials:
\[\mathcal{H}_k = \bigoplus_{p+q=k}\mathcal{H}_{p,q}\]
with eigenvalue $k(k+2n-2)$.
Let $\Delta$ denote the Laplace-Beltrami operator on the sphere quotient $G \backslash S^{2n-1}$. Recall that the eigenvalue of $\Box_b$ on $\mathcal{H}_{p,q}$ is $2q(p+n-1)$. Evidently $\Delta$ and $\Box_b$ are simultaneously diagonalizable. We choose some basis $f_l$ for $(\ker \Box_b)^\perp$, with each $f_l$ an eigenvector of $\Box_b$ and $\Delta$ with eigenvalues $\mu_l$ for $\Delta$, and $\lambda_l$ for $\Box_b$. 
\begin{prop}
    There exists a constant $C$ such that $\forall f \in L^2(G \backslash S^{2n-1})$, $\|\mathcal{G} f\|_{H^s(G \backslash S^{2n-1})}^2 \leq C\|f\|^2_{L^2(G \backslash S^{2n-1})}$ if and only if $\left\{ \frac{(1 + \mu_l)^{\frac{s}{2}}}{\lambda_l} \right\}$ is bounded.
\end{prop}
See \cite{REU2019Sobolev} for details. In \cite{REU2019Sobolev} the authors show that $\mathcal{G}$ gives one weak derivative on spheres, and not more. This should \textit{a priori} imply that for sphere quotients, $\mathcal{G}$ gives exactly one weak derivative. Since one can easily find compact operators that give extra weak derivatives on some sphere quotients, one should be more careful. 
\begin{thm}\label{thm:1 wk derivative}
    There is a constant $C$ such that $\forall f \in H^t(G \backslash S^{2n-1})$, $\|\mathcal{G} f\|_{H^{s+t}(G \backslash S^{2n-1})} \leq C\|f\|_{H^t(G \backslash S^{2n-1})}$ if and only if $s \leq 1$.
\end{thm}
\begin{proof}
In order to show that we get exactly one weak derivative, we will use the generating function associated to a sphere quotient. We have
\[F_G(z,w) = \sum_{p,q}z^pw^q\dim \mathcal{H}_{p,q}^G = \frac 1{|G|}\sum_{g \in G}\frac{1-zw}{\det(z-g)\det(w-\overline{g})} \]
(see \cite{GLR}). Let $e(G)$ denote the exponent of $G$. We observe that
\[F_G(z,w) \frac{(z^{e(G)}-1)^n(w^{e(G)}-1)^n}{1-zw}\]
is a polynomial of degree at most $n(e(G)-1)$ in both variables, and denote it 
\[P_G(z,w) = \sum_{a,b=0}^{n(e(G)-1)}c_{a,b}z^aw^b.\] 
Note this is almost identical to \cite{GLR}, as their formula works for a general group, not just a cyclic one. We replaced the order of $G$ with $e(G)$ as this is enough. In fact, their computation of the coefficients of $P_G$ also works, with a scaling by $\frac {1}{|G|}$ as our definition of $P_G$ is slightly different. We compute
\[\dim \mathcal{H}^G_{0,m e(G)} = \sum_{j=0}^{n-1}\binom{m-j+n-1}{n-1}c_{0,j e(G)}.\]
Note that this formula is valid for all $m \geq 0$, not just $m \geq n$. The right hand side is a polynomial in $m$. We observe $\dim \mathcal{H}^G_{0,0} = 1$, which shows that it is not identically zero. So $\dim \mathcal{H}^G_{0,m e(G)}$ can only be zero for finitely many values of $m$. In particular, $\exists M > 0$ such that for all $m \geq M$,
\[\dim \mathcal{H}^G_{0,m e(G)} \geq 1.\]
The eigenvalue of $\mathcal{H}^G_{0,m e(G)}$ for $\Delta$ is $me(G)(me(G)+2n-2)$, and for $\Box_b$, it is $2me(G)(n-1)$. We see that
\[\frac{\sqrt{1+me(G)(me(G)+2n-2)}^s}{2me(G)(n-1)}\]
is bounded on the set $\{m \in \Z, m \geq M\}$ if and only if $s\leq 1$. Since $H^s(G \backslash S^{2n-1}) \hookrightarrow H^s(S^{2n-1})$, with the norm scaled by $\sqrt{|G|}$, and $\mathcal{G}$ commutes with this injection, evidently $\mathcal{G}$ gains one derivative on $G \backslash S^{2n-1}$ by \cite{REU2019Sobolev}.
\end{proof}
Let 
\[C_{p,q} = \frac{\sqrt{1 +(p+q)(p+q+2n-2)}}{4q(p+n-1)}.\]
Set
\[C_G = \sup_l \frac{\sqrt{1+\mu_l}}{\lambda_l} = \sup_{(p,q), H_{p,q}^G \neq \{0\}}C_{p,q}.\]
Evidently $C_G$ is the optimal constant for \Cref{thm:1 wk derivative}.
It follows by \cite{REU2019Sobolev} that
\[C_G \leq \begin{cases}
    1 & n=2\\
    \frac 12 \frac{\sqrt{n(n-2)}}{(n-1)^2} & n \geq 3
\end{cases}.\]
\begin{thm}
    We have that
    \[C_G > \frac 1{2(n-1)}.\]
    If $G$ is cyclic, then 
    \[C_G \geq \frac{\sqrt{1+4n}}{2n}.\]
\end{thm}
\begin{proof}
    The first inequality follows by the proof of \Cref{thm:1 wk derivative}. Taking the infimum
    \[\inf_{m > 0} \frac{\sqrt{1+me(G)(me(G)+2n-2)}}{2me(G)(n-1)} = \frac 1{2(n-1)}.\]
    We also note that for every $m$, $C_{me(G),0}$ is strictly greater than the infimum. Therefore $C_G$ is strictly greater than the infimum. For the second inequality, note that if $G$ is cyclic, $H_{k,k}^G \neq \{0\}$, as the element $\partial^k_{\overline{z}_n}\partial^k_{z_n}|z|^{2-2n} \in H_{k,k}$ is invariant under diagonal rotations, so is an element of $H_{k,k}^G$. Then we see that $C_{k,k} \leq C_G$ for all $k$. We have
\[C_{k,k}^2 = \frac{1+2k(2k+2n-2)}{4k^2(k+n-1)^2} = \frac{1+2y(k)}{y(k)^2}.\]
Where $y(k) = 2k(k+n-1)$. Set
\[h(y) = \frac{1+2y}{y^2}.\]
Note
\[h'(y) = \frac{2y^2-2y(1+2y)}{y^4}=0\]
when $y = -1$. Hence $h(y)$ is strictly decreasing for $y > 0$. It follows that the maximum constant occurs for the minimal value of $y(k)$. Since $y(k)$ is increasing, this is attained when $k=1$. Hence the maximum value is
\[C_{1,1} = \frac {\sqrt{1+4n}}{2n}.\]
So $C_G \geq \frac {\sqrt{1+4n}}{2n}$.
\end{proof}

\bibliographystyle{plain}
\bibliography{References}

\end{document}